\documentclass[a4paper,12pt]{article}
\usepackage{amssymb}
\usepackage{amsfonts}
\usepackage{amsmath}
\usepackage{geometry}
\usepackage{hyperref}
\usepackage{graphicx} 

\input{tcilatex}
\begin{document}

\title{Analysis of error propagation in the RK3GL2 method}
\author{J. S. C. Prentice \\
Faculty of Information Technology\\
Eduvos\\
Bedfordview, South africa}
\maketitle

\begin{abstract}
The RK3GL2 method is a numerical method for solving initial value problems
in ordinary differential equations, and is a hybrid of a third-order
Runge-Kutta method and two-point Gauss-Legendre quadrature. In this paper we
present an analytical study of the propagation of local errors in this
method, and show that the global order of RK3GL2 is expected to be four.
\end{abstract}

\section{Introduction}

Runge-Kutta (RK) methods are widely used methods for the numerical solution
of initial value problems in ordinary differential equations \cite{hairer}%
\cite{butcher}. Of interest in such methods is the propagation of
approximation error, and the cumulative effect of this propagation. In an RK
method, the accumulation of $O\left( h^{r+1}\right) $ local errors results
in a global error of $O\left( h^{r}\right) ,$ where $h$ is the stepsize. In
other words, the global order of an RK method is one less than its local
order. We have developed a method \cite{prentice 1}, designated RK$r$GL$m$,
which is a hybrid of an RK method of global order $r$, and $m$-point
Gauss-Legendre (GL) quadrature. This method has the property that if the
underlying RK method is $O\left( h^{r+1}\right) $ in its local error then,
with a suitable choice for $m$, the associated RK$r$GL$m$ method is $O\left(
h^{r+1}\right) $ in its \textit{global} error. Recent work has required us
to investigate the nature of local error propagation in the RK3GL2 method.
In this paper we describe such error propagation\ in the RK3GL2 method, and
we show the mechanism by which the global error of RK3GL2 achieves fourth
order. Our study is purely analytical.

\section{Terminology and relevant concepts}

We briefly present notation, terminology and concepts relevant to this study.

\subsection{One-step methods}

We denote an explicit RK method for solving%
\begin{equation*}
y^{\prime }=f\left( x,y\right) \text{ \ \ \ \ \ \ }y\left( x_{0}\right)
=y_{0}\text{ \ \ \ \ \ }a\leqslant x\leqslant b
\end{equation*}%
by%
\begin{equation*}
w_{i+1}=w_{i}+h_{i}F\left( x_{i},w_{i}\right)
\end{equation*}%
where $h_{i}\equiv x_{i+1}-x_{i}$ is a stepsize, $w_{i}$ denotes the
numerical approximation to $y\left( x_{i}\right) $ and $F\left( x,y\right) $
is a function associated with the particular RK method.

For example, the tableau

\renewcommand{\arraystretch}{1.5}%
\begin{equation}
\begin{tabular}{c|ccc}
$0$ &  &  &  \\ 
$\frac{1}{2}$ & $\frac{1}{2}$ &  &  \\ 
$\frac{3}{4}$ & $0$ & $\frac{3}{4}$ &  \\ \hline
& $\frac{2}{9}$ & $\frac{3}{9}$ & $\frac{4}{9}$%
\end{tabular}
\label{tableau}
\end{equation}%
\renewcommand{\arraystretch}{1}has%
\begin{equation*}
hF\left( x,y\right) =\frac{2}{9}k_{1}+\frac{3}{9}k_{2}+\frac{4}{9}k_{3}
\end{equation*}%
and corresponds to a 3rd-order RK method, and is the RK method that we have
used in RK3GL2 in other work. Note that an explicit Runge-Kutta method of
order $r$\ has global error $O\left( h^{r}\right) $ and\ local error $%
O\left( h^{r+1}\right) $\ \cite{butcher}. We use $h$ here as a generic
symbol for the stepsize. We denote the 3rd-order method described here by
RK3. Nevertheless, we note that \emph{any} 3rd-order explicit RK method
could be used in RK3GL2, and we are not restricted to the one presented here.

\subsection{Local and global errors}

We define the global error in a numerical solution at $x_{i}$ by

\begin{equation*}
\Delta _{i}\equiv w_{i}-y_{i},
\end{equation*}%
and the local error at $x_{i}$ by 
\begin{equation}
\varepsilon _{i+1}\equiv \left[ y_{i}+h_{i}F\left( x_{i},y_{i}\right) \right]
-y_{i+1}  \label{eps local}
\end{equation}%
In the above, $y_{i}$ denotes the true solution $y\left( x_{i}\right) .$
Note the use of the exact value $y_{i}$ in the bracketed term in (\ref{eps
local}).

\subsection{Error propagation in a one-step method}

We describe a known result that is useful in our later discussion. We have%
\begin{align*}
w_{1}& =y_{0}+h_{0}F\left( x_{0},y_{0}\right) \\
\Rightarrow \Delta _{1}& =y_{0}+h_{0}F\left( x_{0},y_{0}\right) -y_{1} \\
& =\varepsilon _{1},
\end{align*}%
and%
\begin{align*}
w_{2}& =w_{1}+h_{1}F\left( x_{1},w_{1}\right) \\
\Rightarrow y_{2}+\Delta _{2}& =\left[ y_{1}+\Delta _{1}\right]
+h_{1}F\left( x_{1},y_{1}+\Delta _{1}\right) \\
& =\left[ y_{1}+\Delta _{1}\right] +h_{1}F\left( x_{1},y_{1}\right)
+h_{1}\Delta _{1}F_{y}\left( x_{1},\xi _{1}\right) \\
\Rightarrow \Delta _{2}& =\left[ y_{1}+h_{1}F\left( x_{1},y_{1}\right) -y_{2}%
\right] +\Delta _{1}\left[ 1+h_{1}F_{y}\left( x_{1},\xi _{1}\right) \right]
\\
& =\varepsilon _{2}+\alpha _{1}\varepsilon _{1}.
\end{align*}%
Furthermore,%
\begin{align*}
\Delta _{3}& =\varepsilon _{3}+\alpha _{2}\Delta _{2}=\varepsilon
_{3}+\alpha _{2}\left( \varepsilon _{2}+\alpha _{1}\varepsilon _{1}\right) \\
& =\varepsilon _{3}+\alpha _{2}\varepsilon _{2}+\alpha _{2}\alpha
_{1}\varepsilon _{1}
\end{align*}%
and%
\begin{align*}
\Delta _{4}& =\varepsilon _{4}+\alpha _{3}\Delta _{3} \\
& =\varepsilon _{4}+\alpha _{3}\varepsilon _{3}+\alpha _{3}\alpha
_{2}\varepsilon _{2}+\alpha _{3}\alpha _{2}\alpha _{1}\varepsilon _{1}
\end{align*}%
where%
\begin{align*}
\alpha _{k}& =1+h_{k}F_{y}\left( x_{k},\xi _{k}\right) \\
\xi _{k}& \in \left( y_{k},y_{k}+\Delta _{k}\right) .
\end{align*}%
We present $F_{y}\left( x,y\right) =\partial F\left( x,y\right) /\partial y$
for the RK method (\ref{tableau}) in the Appendix.

In general%
\begin{equation*}
\Delta _{n}=\dsum\limits_{j=1}^{n}\left[ \left( \frac{1}{\alpha _{n}}\right)
\left( \dprod\limits_{k=j}^{n}\alpha _{k}\right) \right] \varepsilon _{j}.
\end{equation*}%
If $\left\vert h_{k}F_{y}\left( x_{k},\xi _{k}\right) \right\vert $\ is
small then $\alpha _{k}\approx 1,$\ and so%
\begin{equation*}
\Delta _{n}\approx \dsum\limits_{j=1}^{n}\varepsilon _{j}
\end{equation*}%
but this is generally not expected to be the case, particularly if $%
F_{y}\left( x_{k},\xi _{k}\right) $\ is large. Furthermore, if the $\alpha $%
's have magnitude larger than unity, then the term in $\varepsilon _{1}$\
could make the most significant contribution to the global error.

The global error $\Delta _{n}$ is the accumulation of these local errors, as
in%
\begin{equation*}
\Delta _{n}=\dsum\limits_{j=1}^{n}\underset{\varepsilon _{j}}{\underbrace{%
\beta _{j}h^{r+1}}}=\left( \frac{1}{n}\dsum\limits_{j=1}^{n}\beta
_{j}\right) \left( nh\right) h^{r}=\overline{\beta }\left( b-a\right) h^{r}
\end{equation*}%
where $\beta _{j}$ and $\overline{\beta }$ have been implicitly defined, and
we have used $nh=b-a$ (so that, in this expression, $h$ is the average
separation of the nodes $x_{j}).$ Note that $\Delta _{n}$ is $O\left(
h^{r}\right) .$

\subsection{Gauss-Legendre quadrature}

Gauss-Legendre (GL) quadrature on $[u,v]$\ with $2$\ nodes is given by \cite%
{kincaid}%
\begin{equation*}
\dint\limits_{u}^{v}f\left( x,y\right)
dx=h\dsum\limits_{i=1}^{2}C_{i}f\left( x_{i},y_{i}\right) +O\left(
h^{5}\right)
\end{equation*}%
where the nodes $x_{i}$\ are the roots of the $2$nd degree Legendre
polynomial on $[u,v].$ Here, $h$ is the average separation of the nodes on $%
[u,v]$, a notation we will adopt from now on, and the $C_{i}$\ are
appropriate weights. For GL2, the roots of the 2nd degree Legendre
polynomial on $[-1,1]$ are 
\begin{equation*}
\widetilde{x}_{1}=-\frac{\sqrt{3}}{3},\text{ \ }\widetilde{x}_{2}=\frac{%
\sqrt{3}}{3}
\end{equation*}%
and are mapped to corresponding nodes $x_{i}$ on $[u,v]$ via%
\begin{equation*}
x_{i}=\frac{1}{2}\left[ \left( v-u\right) \widetilde{x}_{i}+u+v\right] .
\end{equation*}%
Also, the average node separation on $[-1,1]$ is $2/3$, and so $h$ on $[u,v]$
is given by%
\begin{equation*}
h=\frac{2}{3}\left( \frac{v-u}{2}\right) =\frac{v-u}{3},
\end{equation*}%
while the weights%
\begin{equation*}
C_{1}=\frac{3}{2},\text{ \ }C_{2}=\frac{3}{2}
\end{equation*}%
are constants on any interval of integration.

\begin{figure}
    \centering
    \includegraphics[width=1\linewidth]{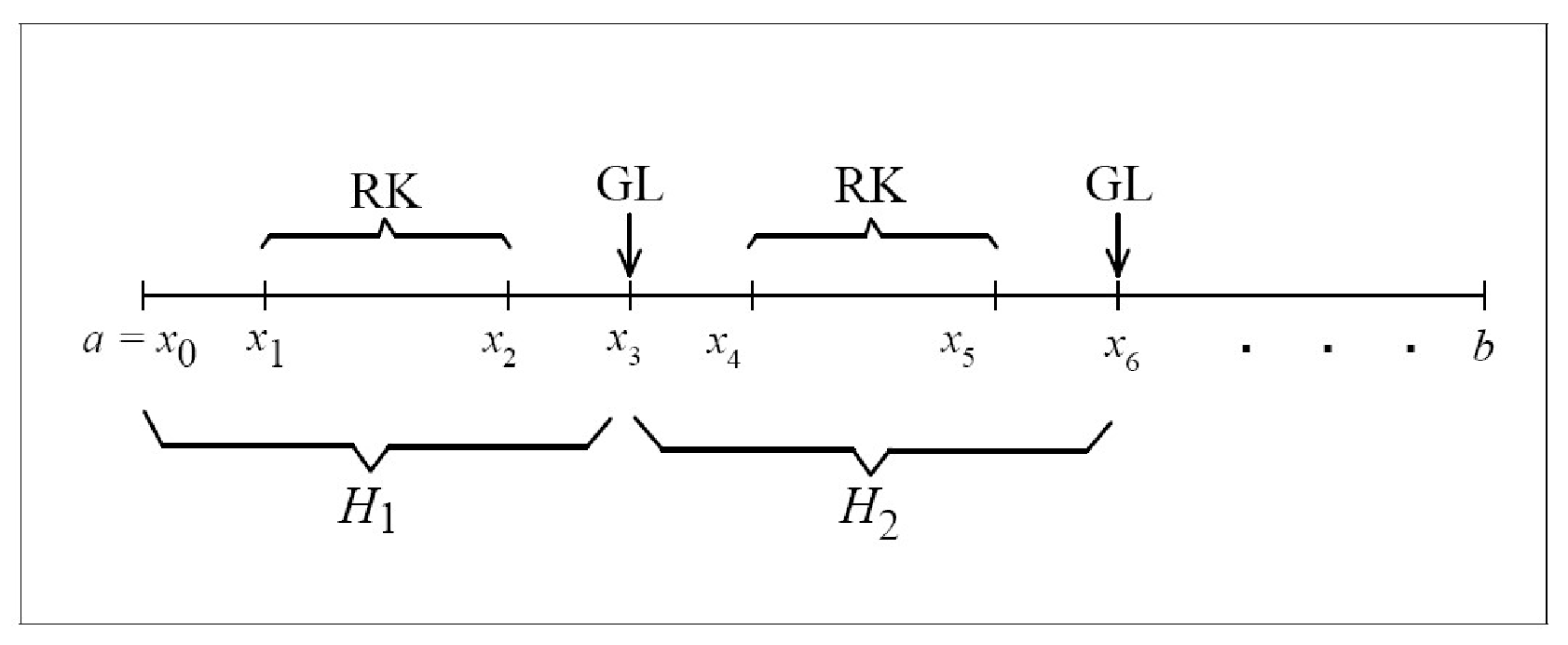}
    \caption{RK3GL2 algorithm for the first two subintervals $H_{1}$ and $H_{2}$ on $[a,b].$}
    \label{figure 1}
\end{figure}

\subsection{The RK3GL2 algorithm}

We briefly describe the general RK3GL2 algorithm \cite{prentice 1} on the
interval $[a,b]$, with reference to Figure 1.

Subdivide $[a,b]$ into $N$\
subintervals $H_{i}.$ At the RK nodes we use RK3$:$%
\begin{equation*}
w_{i+1}=w_{i}+h_{i}F\left( x_{i},w_{i}\right) 
\end{equation*}%
At the GL nodes we use $2$-point GL quadrature:%
\begin{align*}
w_{3}& =w_{0}+h\dsum\limits_{j=1}^{2}C_{j}f\left( x_{j},w_{j}\right) , \\
w_{6}& =w_{3}+h\dsum\limits_{j=4}^{5}C_{j}f\left( x_{j},w_{j}\right) ,
\end{align*}%
and so on. The GL component is motivated by%
\begin{align*}
\dint\limits_{x_{0}}^{x_{3}}f\left( x,y\right) dx& =y_{3}-y_{0}\approx
h\dsum\limits_{j=1}^{2}C_{j}f\left( x_{j},y_{j}\right) \Rightarrow
y_{3}\approx y_{0}+h\dsum\limits_{j=1}^{2}C_{j}f\left( x_{j},y_{j}\right) ,
\\
\dint\limits_{x_{3}}^{x_{6}}f\left( x,y\right) dx& =y_{6}-y_{3}\approx
h\dsum\limits_{j=4}^{5}C_{j}f\left( x_{j},y_{j}\right) \Rightarrow
y_{6}\approx y_{3}+h\dsum\limits_{j=4}^{5}C_{j}f\left( x_{j},y_{j}\right) ,
\end{align*}%
and so on. Of course, RK3GL2 is RK$r$GL$m$ with $r=3,m=2$. The RK3GL2
algorithm has been shown to be consistent, convergent and zero-stable \cite%
{prentice 1}.

\subsection{Local error at the GL nodes}

The local error at the GL nodes is defined in a similar way to that for a
one-step method: for example,%
\begin{align*}
\dint\limits_{x_{0}}^{x_{3}}f\left( x,y\right) dx&
=y_{3}-y_{0}=h\dsum\limits_{j=1}^{2}C_{j}f\left( x_{j},y_{j}\right) +O\left(
h^{5}\right) \\
\Rightarrow \varepsilon _{3}& =\underset{\text{exact values of }y\left(
x\right) }{\underbrace{\left[ y_{0}+h\dsum\limits_{j=1}^{2}C_{j}f\left(
x_{j},y_{j}\right) \right] }}-y_{3}=O\left( h^{5}\right) .
\end{align*}%
In RK3GL2, then, the local error at the GL nodes is $O\left( h^{5}\right) .$

\section{Error propagation in the RKGL method}

\subsection{Error propagation}

For RK3GL2 we have

\begin{align}
\Delta _{1}& =\varepsilon _{1},\text{ \ \ }\Delta _{2}=\varepsilon
_{2}+\alpha _{1}\varepsilon _{1}  \label{delta1 delta2 delta3} \\
w_{3}=y_{3}+\Delta _{3}& =y_{0}+h\dsum\limits_{j=1}^{2}C_{j}f\left(
x_{j},w_{j}\right) =y_{0}+h\dsum\limits_{j=1}^{2}C_{j}f\left(
x_{j},y_{j}+\Delta _{j}\right)  \notag \\
& =y_{0}+h\dsum\limits_{j=1}^{2}C_{j}f\left( x_{j},y_{j}\right)
+h\dsum\limits_{j=1}^{2}C_{j}f_{y}\left( x_{j},\zeta _{j}\right) \Delta _{j}
\notag \\
& =y_{0}+h\dsum\limits_{j=1}^{2}C_{j}f\left( x_{j},y_{j}\right)
+h\dsum\limits_{j=1}^{2}\gamma _{j}\varepsilon _{j}  \notag \\
& =y_{0}+h\dsum\limits_{j=1}^{2}C_{j}f\left( x_{j},y_{j}\right) +A_{1,2}h 
\notag \\
\Rightarrow \Delta _{3}& =\left[ y_{0}+h\dsum\limits_{j=1}^{2}C_{j}f\left(
x_{j},y_{j}\right) -y_{3}\right] +A_{1,2}h  \notag \\
& =\varepsilon _{3}+A_{1,2}h  \notag
\end{align}%
where $\zeta _{j}\in \left( y_{j},y_{j}+\Delta _{j}\right) ,$ and $\gamma
_{j}$ and $A_{1,2}$ have been implicitly defined. Indeed, $\gamma _{j}$
arises from the fact that $\Delta _{1}$ and $\Delta _{2}$ can be written in
terms of $\varepsilon _{1}$ and $\varepsilon _{2}.$ Also,%
\begin{equation}
\Delta _{4}=\varepsilon _{4}+\alpha _{3}\Delta _{3},\text{ \ \ }\Delta
_{5}=\varepsilon _{5}+\alpha _{4}\varepsilon _{4}+\alpha _{4}\alpha
_{3}\Delta _{3}  \label{delta5 delta6}
\end{equation}%
so that%
\begin{align*}
w_{6}=y_{6}+\Delta _{6}& =w_{3}+h\dsum\limits_{j=4}^{5}C_{j}f\left(
x_{j},w_{j}\right) =y_{3}+\Delta _{3}+h\dsum\limits_{j=4}^{5}C_{j}f\left(
x_{j},y_{j}+\Delta _{j}\right) \\
& =y_{3}+\Delta _{3}+h\dsum\limits_{j=4}^{5}\left[ C_{j}f\left(
x_{j},y_{j}\right) +C_{j}f_{y}\left( x_{j},\zeta _{j}\right) \Delta _{j}%
\right] \\
& =y_{3}+h\dsum\limits_{j=4}^{5}C_{j}f\left( x_{j},y_{j}\right)
+A_{4,5}h+B_{6}\Delta _{3}h+\Delta _{3} \\
\Rightarrow \Delta _{6}& =\left[ y_{3}+h\dsum\limits_{j=4}^{5}C_{j}f\left(
x_{j},y_{j}\right) -y_{6}\right] +A_{4,5}h+B_{6}\Delta _{3}h+\varepsilon
_{3}+A_{1,2}h \\
& =\left( \varepsilon _{6}+\varepsilon _{3}\right) +\left(
A_{4,5}h+A_{1,2}h\right) +B_{6}\Delta _{3}h,
\end{align*}%
where $A_{4,5}$ and $B_{6}$ are defined below in (\ref{A4N-3,4N-1}) and (\ref%
{B8}). In fact, $B_{6}$ arises from the fact that $\Delta _{4}$ and $\Delta
_{5}$ can be written in terms of $\Delta _{3}.$Furthermore,%
\begin{equation*}
\Delta _{9}=\left( \varepsilon _{9}+\varepsilon _{6}+\varepsilon _{3}\right)
+\left( A_{7,8}h+A_{4,5}h+A_{1,2}h\right) +B_{9}\Delta _{6}h+B_{6}\Delta
_{3}h.
\end{equation*}%
and%
\begin{align*}
\Delta _{12}& =\left( \varepsilon _{12}+\varepsilon _{9}+\varepsilon
_{6}+\varepsilon _{3}\right) +\left(
A_{10,11}h+A_{7,8}h+A_{4,5}h+A_{1,2}h\right) \\
& +B_{12}\Delta _{9}h+B_{9}\Delta _{6}h+B_{6}\Delta _{3}h.
\end{align*}%
In general,%
\begin{align*}
\Delta _{3N}=& \left( \varepsilon _{3N}+\cdots +\varepsilon _{3}\right)
+\left( A_{3N-2,3N-1}h+\cdots +A_{1,2}h\right) \\
& +\left( B_{3N}\Delta _{3\left( N-1\right) }h+B_{3\left( N-1\right) }\Delta
_{3\left( N-2\right) }h+\cdots +B_{6}\Delta _{3}h\right) .
\end{align*}%
In the above, $A_{1,2},A_{4,5},A_{7,8},A_{3N-2,3N-1},B_{6},B_{9}$ and $%
B_{12} $\ are appropriate coefficients (defined below), and $N$\ is the
total number of subintervals into which $[a,b]$\ has been subdivided. We
list below some relevant terms in detail.%
\begin{align}
A_{3N-2,3N-1}& =\dsum\limits_{j=3N-2}^{3N-1}\gamma _{j}\varepsilon _{j}.
\label{A4N-3,4N-1} \\
\gamma _{3N-2}& =\left( 1+\alpha _{3N-2}\right) C_{3N-2}f_{y}\left(
x_{3N-2},\zeta _{3N-2}\right) .  \notag \\
\gamma _{3N-1}& =C_{3N-1}f_{y}\left( x_{3N-1},\zeta _{3N-1}\right) .  \notag
\end{align}%
\begin{align}
B_{6}& =C_{4}f_{y}\left( x_{4},\zeta _{4}\right) \alpha
_{3}+C_{5}f_{y}\left( x_{5},\zeta _{5}\right) \alpha _{3}\alpha _{4}.
\label{B8} \\
B_{9}& =C_{7}f_{y}\left( x_{7},\zeta _{7}\right) \alpha
_{6}+C_{8}f_{y}\left( x_{8},\zeta _{8}\right) \alpha _{6}\alpha _{7}.  \notag
\\
B_{12}& =C_{10}f_{y}\left( x_{10},\zeta _{10}\right) \alpha
_{9}+C_{11}f_{y}\left( x_{11},\zeta _{11}\right) \alpha _{9}\alpha _{10}. 
\notag
\end{align}%
In general,%
\begin{equation*}
B_{3N}=C_{3N-2}f_{y}\left( x_{3N-2},\zeta _{3N-2}\right) \alpha
_{3N-3}+C_{3N-1}f_{y}\left( x_{3N-1},\zeta _{3N-1}\right) \alpha
_{3N-3}\alpha _{3N-2}.
\end{equation*}%
For completeness, we could include a term of the form%
\begin{equation*}
B_{3}\Delta _{0}h
\end{equation*}%
in the expression for $\Delta _{3N}$ above, but we assume that $\Delta
_{0}=0 $ so that such a term is not necessary here.

As for writing the global error in terms of the local errors consider, for
example, $\Delta _{12}:$

\begin{align*}
\Delta _{12}=& \left( \varepsilon _{12}+\varepsilon _{9}+\varepsilon
_{6}+\varepsilon _{3}\right) +\left(
A_{10,11}+A_{7,8}+A_{4,5}+A_{1,2}\right) h+\left( B_{12}\Delta
_{9}+B_{9}\Delta _{6}+B_{6}\Delta _{3}\right) h \\
=& \varepsilon _{12}+\varepsilon _{9}+\varepsilon _{6}+\varepsilon
_{3}+\dsum\limits_{j=10}^{11}\gamma _{j}h\varepsilon
_{j}+\dsum\limits_{j=7}^{8}\gamma _{j}h\varepsilon
_{j}+\dsum\limits_{j=4}^{5}\gamma _{j}h\varepsilon
_{j}+\dsum\limits_{j=1}^{2}\gamma _{j}h\varepsilon _{j} \\
& +B_{12}h\varepsilon _{9}+B_{12}h\varepsilon _{6}+B_{12}h\varepsilon
_{3}+\dsum\limits_{j=7}^{8}B_{12}\gamma _{j}h^{2}\varepsilon
_{j}+\dsum\limits_{j=4}^{5}B_{12}\gamma _{j}h^{2}\varepsilon
_{j}+\dsum\limits_{j=1}^{2}B_{12}\gamma _{j}h^{2}\varepsilon _{j} \\
& +B_{12}B_{9}h^{2}\varepsilon _{6}+B_{12}B_{9}h^{2}\varepsilon
_{3}+\dsum\limits_{j=4}^{5}B_{12}B_{9}\gamma _{j}h^{3}\varepsilon
_{j}+\dsum\limits_{j=1}^{2}B_{12}B_{9}\gamma _{j}h^{3}\varepsilon _{j} \\
& +B_{12}B_{9}B_{6}h^{3}\varepsilon
_{3}+\dsum\limits_{j=1}^{2}B_{12}B_{9}B_{6}\gamma _{j}h^{4}\varepsilon
_{j}+B_{9}h\varepsilon _{6}+B_{9}h\varepsilon
_{3}+\dsum\limits_{j=4}^{5}B_{9}\gamma _{j}h^{2}\varepsilon _{j} \\
& +\dsum\limits_{j=1}^{2}B_{9}\gamma _{j}h^{2}\varepsilon
_{j}+B_{9}B_{6}h^{2}\varepsilon _{3}+\dsum\limits_{j=1}^{2}B_{9}B_{6}\gamma
_{j}h^{3}\varepsilon _{j}+B_{6}h\varepsilon
_{3}+\dsum\limits_{j=1}^{2}B_{6}\gamma _{j}h^{2}\varepsilon _{j}.
\end{align*}%
Using the above expressions, we have, in terms of the local errors $%
\varepsilon _{i},$%
\begin{equation*}
\Delta _{12}=\sum_{i=1}^{12}G_{i}\varepsilon _{i}
\end{equation*}%
where%
\begin{align}
G_{1}& =\gamma _{1}h+B_{12}\gamma _{1}h^{2}+B_{9}\gamma
_{1}h^{2}+B_{6}\gamma _{1}h^{2}+B_{12}B_{9}\gamma _{1}h^{3}+B_{9}B_{6}\gamma
_{1}h^{3}+B_{12}B_{9}B_{6}\gamma _{1}h^{4}  \notag \\
G_{2}& =\gamma _{2}h+B_{12}\gamma _{2}h^{2}+B_{9}\gamma
_{2}h^{2}+B_{6}\gamma _{2}h^{2}+B_{12}B_{9}\gamma _{2}h^{3}+B_{9}B_{6}\gamma
_{2}h^{3}+B_{12}B_{9}B_{6}\gamma _{2}h^{4}  \notag \\
G_{3}&
=1+B_{12}h+B_{9}h+B_{6}h+B_{12}B_{9}h^{2}+B_{9}B_{6}h^{2}+B_{12}B_{9}B_{6}h^{3}
\notag \\
G_{4}& =\gamma _{4}h+B_{12}\gamma _{4}h^{2}+B_{9}\gamma
_{4}h^{2}+B_{12}B_{9}\gamma _{4}h^{3}  \notag \\
G_{5}& =\gamma _{5}h+B_{12}\gamma _{5}h^{2}+B_{9}\gamma
_{5}h^{2}+B_{12}B_{9}\gamma _{5}h^{3}  \notag \\
G_{6}& =1+B_{12}h+B_{9}h+B_{12}B_{9}h^{2}  \label{g6} \\
G_{7}& =\gamma _{7}h+B_{12}\gamma _{7}h^{2}  \notag \\
G_{8}& =\gamma _{8}h+B_{12}\gamma _{8}h^{2}  \notag \\
G_{9}& =1+B_{12}h  \notag \\
\text{ }G_{10}& =\gamma _{10}h  \notag \\
G_{11}& =\gamma _{11}h  \notag \\
G_{12}& =1  \notag
\end{align}%
The global error at the RK nodes is understood with reference to section
2.3, and equations (\ref{delta1 delta2 delta3}) and (\ref{delta5 delta6}).

\subsection{Error accumulation}

We have%
\begin{align}
\Delta _{3N}=& \left( \varepsilon _{3N}+\cdots +\varepsilon _{3}\right)
+\left( A_{3N-2,3N-1}h+\cdots +A_{1,2}h\right)  \notag \\
& +\left( B_{3N}\Delta _{3\left( N-1\right) }h+\cdots +B_{6}\Delta
_{3}h\right)  \label{qwqw} \\
=& NM_{1}h^{5}+NM_{2}h^{6}  \notag \\
=& \left( \frac{M_{1}}{3}\right) h^{4}\left[ 3Nh\right] +\left( \frac{M_{2}}{%
3}\right) h^{5}\left[ 3Nh\right]  \notag \\
=& \left[ \frac{M_{1}\left( b-a\right) }{3}\right] h^{4}+\left[ \frac{%
M_{2}\left( b-a\right) }{3}\right] h^{5}  \notag \\
=& O\left( h^{4}\right) ,  \notag
\end{align}%
where $M_{1}$ and $M_{2}$ are appropriate coefficients. This demonstrates
the $O\left( h^{4}\right) $ character of the global error in RK3GL2.

\section{Comments}

The mechanism for the $O\left( h^{4}\right) $ global error in RK3GL2 is
shown in the first two terms on the RHS of (\ref{qwqw}). The first of these
is the sum of the GL local errors which, through a suitable choice of $m$,
is $\propto h^{5}$. The second term is a linear combination of the RK local
errors, multiplied by a factor $h.$ Since each RK local error is $\propto
h^{5}$, this term is $\propto h^{6}$. The effect of the GL component, then,
is to increase the order of the accumulated RK local errors by one. We refer
to this as a \emph{quenching}\ effect that occurs at the GL nodes, and it
serves to damp the accumulation of the RK local errors. The third term in (%
\ref{qwqw}) contains terms $\propto h^{6}$ and higher, as shown, for
example, in the expansion of $\Delta _{12}$ via $G_{1},...,G_{12}$ in (\ref%
{g6}).

\begin{flushleft}
\textbf{{\Large {Appendix }}}

\ \ \ \ \ \ \ \ \ \ \ \ 
\end{flushleft}

As a matter of interest, with%
\begin{align*}
k_{1}& =hf\left( x,y\right) \\
k_{2}& =hf\left( x+\frac{h}{2},y+\frac{k_{1}}{2}\right) =hf\left( x+\frac{h}{%
2},y+\frac{hf\left( x,y\right) }{2}\right) \\
k_{3}& =hf\left( x+\frac{3h}{4},y+\frac{3k_{2}}{4}\right) =hf\left( x+\frac{%
3h}{4},y+\frac{3}{4}hf\left( x+\frac{h}{2},y+\frac{hf\left( x,y\right) }{2}%
\right) \right) ,
\end{align*}%
we find%
\begin{align*}
F_{y}\left( x,y\right) & =\frac{2hf_{y}\left( x,y\right) }{9}+\frac{%
3hf_{y}\left( x+\frac{h}{2},y+\frac{k_{1}}{2}\right) }{9}\left( 1+\frac{%
hf_{y}\left( x,y\right) }{2}\right) \\
& +\frac{4hf_{y}\left( x+\frac{3h}{4},y+\frac{3k_{2}}{4}\right) }{9}\left( 1+%
\frac{3hf_{y}\left( x+\frac{3h}{4},y+\frac{k_{1}}{2}\right) }{4}\left( 1+%
\frac{hf_{y}\left( x,y\right) }{2}\right) \right) .
\end{align*}

\end{document}